\documentclass[12pt]{article}
\usepackage{cite,epsfig}
\newcommand{\be}{\begin{equation}}
\newcommand{\ee}{\end{equation}}
\newcommand{\ben}{\begin{eqnarray}}
\newcommand{\een}{\end{eqnarray}}
\newcommand{\un}{\underline}

\newcommand{\arexn}{\imath\frac{ n \pi t }{T}}

\newcommand{\arexna}{\imath\frac{a n\pi t }{T}}
\newcommand{\arexma}{\imath\frac{a m \pi t}{T}}

\newcommand{\exna}{\frac{1}{\sqrt{2T}}e^{\arexna}}
\newcommand{\exnachi}{\frac{\chi_T(t)}{\sqrt{2T}}e^{\arexna}}

\newcommand{\suni}{\sum^{\infty}_{n=-\infty}}

\newcommand{\intT}{\int_{-T}^{T}}

\newcommand{\vcp}{\vec{c'}}

\begin{document}
\baselineskip = 2\baselineskip

\title{Oversampling and transmission of hidden information}
\author{Jody R. Miotke and  Laura Rebollo-Neira\\
NCRG, Aston University,\\
Birmingham B4 7ET,\\
United Kingdom\\
http:$//$www.ncrg.aston.ac.uk}
\date{}


\maketitle

\begin{abstract}
The lack of uniqueness arising by oversampling 
of Fourier coefficients is shown to provide  
a way of transmitting hidden information. 
A basic encoding/decoding
system, developed on the basis of such a possibility, 
is discussed. 
The system is devised with the double purpose of: 
(a) enabling the transmission of an arbitrary signal and 
(b) allowing for the transmission of a hidden embedded code. 
\end{abstract}


\section{Introduction}
In a previous publication \cite{reco} the oversampling 
problem of signal representation by
Discrete Fourier Transform (DFT)
has been addressed from the Frame Theory \cite{yo,do,ol}
point of view. 
Within this structure the 
signal reconstruction 
appears as a {\it{tight frame}} superposition.  
As remarked in \cite{reco}, the implicit redundancy of
frames is 
the cause of noise reduction in signal reconstruction 
 and also the reason
that the representation is {\it{not unique}}.
In this communication 
we show that the lack of uniqueness can be  used
for dissembling information. \\
We propose a basic 
encoding/decoding 
procedure, which enables the transmission of   hidden 
information when transmitting some arbitrary signal.
Although we shall restrict our  considerations to 
redundancy arising 
by DFT oversampling,  
the proposal for transmitting hidden information could 
be equivalently implemented using any other 
redundant transformation. 
In particular, recently introduced techniques for 
the construction of tight frames 
\cite{c1,c2,c3,c4} 
 should be relevant to development of more 
sophisticated systems based 
on the principles proposed here. Nevertheless,  
we feel it is appropriate to 
present the central ideas of 
this contribution by making use of the 
oversampled DFT. Among a number of   
important reasons leading to this choice we would like 
to point out the following: 
i)The Fast Fourier 
Transform (FFT), fast implementation of DFT, 
is still without  doubt the most 
popular and widely used signal processing tool. 
ii)Oversampling FFT, by  padding with zeros, 
is a common procedure for noise reduction and smoothing, 
hence the importance of stressing that 
it also leaves room for the transmission of hidden information.
iii)All readers, regardless of their 
area of expertise, can be assumed to be familiar with the 
corresponding mathematical background.  
Consequently, we shall introduce the proposed technique within 
the framework arising by oversampling DFT coefficients.
The effect of additive noise in the transmission channel 
is also considered. In order to discriminate the level of
additive zero-mean noise
in which the system can safely operate,
two possible situations are analysed.\\ 
The letter is organised as follows: In Section 2 the lack 
of uniqueness 
inherent to oversampling, as well at its relevance 
in relation  to  
transmission of hidden information, are discussed. A 
basic DFT based encoding/decoding system is proposed in 
Section 3, and some 
examples to illustrate its performance in the presence of 
additive zero-mean 
 noisy are given in Section 4. 
 The conclusions are drawn in Section 5.  
\section{Oversampling and lack of uniqueness}
Let us  represent a  signal $f(t)$, which is defined for 
$t \in [-T,T]$,   through  its 
 Discrete Fourier expansion, i.e.
\be
f(t) = \frac{1}{\sqrt{2T}}
\suni c_n  e^{\arexn}.
\label{one}
\ee
Since for $t \in [-T,T]$ the complex exponentials in 
(\ref{one})
constitute an orthonormal basis,  
the coefficients $c_n$ in (\ref{one}) are obtained 
as: 
\be
c_n= \frac{{1}}{\sqrt{2T}}\intT f(t) e^{-\arexn}\, dt.
\ee
Let us consider now the re-scaling operation: $t \to a t$, with
$a$ a positive  real number less than 1, and construct the 
functions $\exnachi$, with $\chi_T(t)$ 
 deffined as: 
 $\chi_T(t)=1$ if $t \in [-T,T]$ and
zero otherwise. The new functions $\exnachi$ are no 
longer a basis but a {\it{tight frame}}  for the 
space of time limited signal with time-width $T$ (the 
corresponding frame-bound being $a^{-1}$ \cite{reco}).
This has 
a remarkable consequence, namely the coefficients 
$c_n$ of the linear expansion:
\be
f(t) = \frac{{\chi_T(t)}}{\sqrt{2T}} \suni c_n  e^{\arexna}
\label{tf}
\ee
are {\it{not unique}}. There 
exist infinitely
many different sets of coefficients $c_n$ which
can reproduce an identical signal $f$ by the above linear
superposition. A {{\un{particular}} set of
coefficients $c_n$ is obtained as:
\be
c_n =  \frac{{a}}{\sqrt{2T}} \intT f(t) e^{-\arexna} 
\, dt.
\label{mn}
\ee
Out of all possible sets of coefficients, the ones 
given by the above equation constitute 
the coefficients of minimum 2-norm \cite{yo,do}.\\
Let us stress the cause for the coefficients
in the tight frame  expansion not to be unique. 
The reason being that, for $a<1$ with the 
restriction $t \in [-T,T]$ the exponentials
$\exna$ are {\it{not}} linearly independent i.e.,
we can have the situation:
$$\frac{{1}}{\sqrt{2T}} \suni c'_n   e^{\arexna} =0 
\,\,\,\,\, {\mbox{for}}\,\, \,\,\,
\suni |c'_n|^2  \neq 0$$
or, taking inner products both sides  with
$\frac{1}{\sqrt{2T}}{e^{\arexma}}$,
$$\frac{1}{{2T}}\suni c'_n \intT e^{-{\arexma}}  e^{\arexna}
\, dt =0 \,\,\,\,\, {\mbox{for}}\,\, \,\,\,
\suni |c'_n|^2  \neq 0,$$
which can be recast in the fashion:
$$ G \vcp =0 \,\,\,\,\, {\mbox{for}}\,\, \,\,\,
 || \vec{c'} ||^ 2 = \suni |c'_n|^2  \neq 0,$$
with $G$ a matrix of elements:
\be
g_{m,n}=\frac{1}{2T} \intT e^{-{\arexma}} e^{\arexna}
 \, dt= \frac{\sin {a (m-n)\pi}}{a(m-n) \pi} 
\label{gmn}
\ee
and $\vcp$ a
vector of  components $c_n'$.\\
Notice that all vectors $\vcp$ satisfying
$G \vcp =0 $ belong, by definition,
to Null($G$),
the Null
space of $G$. All
such vectors satisfy:
\be
f(t) = \frac{{\chi_T(t)}}{\sqrt{2T}}  \suni c_n   e^{\arexna} + 
 \frac{{\chi_T(t)}}{\sqrt{2T}} \suni c'_n  e^{\arexna}  =  
  \frac{{\chi_T(t)}}{\sqrt{2T}}  \suni  c_n'' e^{\arexna},
\label{six}
\ee
where we have defined  $c''_n  = c_n + c'_n$ with 
$c_n$ as in (\ref{mn}) and $c'_n$ the components 
of an arbitrary 
vector $\vcp \in$ Null($G$). Vectors 
$\vec{c}$ and $\vcp$ will hereafter  
be referred to as {\it {signal coefficients}}  
and {\it {hidden code coefficients}}
respectively. The fact that all coefficients 
$\vec{c''}= \vec{c} + \vcp$ 
reproduce an identical signal as coefficients $\vec{c}$
provides us with the foundation to construct 
an encoding/decoding
scheme for transmitting hidden information. 

\section{The encoding-decoding system}

Let us assume that, in addition to transmitting an arbitrary 
signal $f$, we wish to transmit a hidden code $\vec{h}$ 
consisting of $K$ numbers. 
For practical implementation we give to the 
oversampling parameter $a$ a positive value, less than one, and 
consider that
$G$ is an $M\times M$  matrix of elements
as  given in (\ref{gmn}). 
We select $K$ eigenvectors of $G$
corresponding to the zero eigenvalues, which are assumed to be 
orthonormal, and construct a 
vector $\vec{c'}\in $ Null($G$) as follows: 
\be
\vec{c'}= U  B_s \vec{h}
\label{hi}
\ee
where $U$ is an $M\times K$ matrix, the  columns of which
are the $K$ selected eigenvectors and $B_s$ is a $K \times K$
unitary random matrix. Note: the subindex $s$ indicates  
that the random generator used for constructing the matrix 
is initialized at state  $s$. Such a state is needed to be 
known at the decoding stage. 

\subsection*{Encoding process}
Consider that the signal $f$ to be transmitted 
is given as an $N$-dimension data vector and proceed as 
follows: 

\begin{itemize}

\item 
Compute the signal coefficients $\vec{c}$ as in (\ref{mn}).\\
Note that this calculation 
can be carried out with Fast Fourier Transform (FFT)
by adding $\frac{N (1-a)}{2a}$ zeros at the beginning and 
at the end of the data vector $f$, so as to obtain the 
required vector $\vec{c}$ of dimension $M=\frac{N}{a}$.

\item Compute the hidden code coefficients $\vec{c'}$ as 
prescribed in (\ref{hi}).  

\item Transmit the coefficients 
$\vec{c''}=\vec{c} + \vec{c'}$ to the 
receiver. 

\end{itemize}

{\subsection*{Decoding process}}
\begin{itemize}
\item 
Use the received vector  $\vec{c''}$ 
for recovering the signal $f$.\\
In practice this  can be 
 computed by Inverse Fast Fourier Transform (IFFT)
on the received vector $\vec{c''}$. 

\item 
Use the signal $f$ to compute the 
 signal coefficients $\vec{c}$ as in (\ref{mn}), which can be 
  accomplished by applying  FFT on the signal 
$f$ recovered in the previous step.

\item  
Compute vector $\vec{c'}$ 
through $\vec{c'}= \vec{c''} - \vec{c}$. 

\item  
Recover the hidden  code $\vec{h}$ by noticing that:

(a)For constructing the matrix $U$ one can use all eigenvectors of 
 the matrix $G$ corresponding to eigenvalues less than a previously 
specified tolerance  parameter.
Matrix $U$ is unitary, i.e., 
${U}^{-1}= U^{*}$ and then we have: 
$B_s \vec{h} = U^{*} \vec{c'}$, 
where $U^{*}$ indicates the transpose conjugate 
of matrix $U$. 

(b)The dimension of matrix $B_s$ can be determined from 
the number of  non-zero  components of vector $U^{*}\vec{c'}$. 
Thereby, state $s$ allows the  
reproduction of the random matrix $B_s$. 
Since this is also a unitary matrix  
${B_s}^{-1}= B_s^{*}$. 

Hence the  vector $\vec{h}$ is obtained as:

$$\vec{h} = {B_s^{*}} U^{*} \vec{c'}.$$

\end{itemize}

{\it{Remarks:}} In order to be able to implement the 
above described decoding process the receiver should know:\\
i)The transmitted vector $\vec{c''}$.\\
ii)The oversampling parameter $a$, yet this parameter 
might be actually estimated in some situations
by counting the 
beginning/ending zeros  apparent in the recovered signal.\\
iii)The state $s$ that was used for generating the unitary 
matrix $B_s$. Note: the restriction of unitariness is
imposed in 
order to avoid amplification of errors in the inversion process. 
Such a condition is achieved by orthogonalization of a  
well posed random matrix.\\

With the knowledge of i) ii) and iii) the receiver should be 
able to reproduce both the signal $f$ and the hidden code $\vec{h}$.

\section{Example}
Consider that we wish to transmit the chirp signal of Figure 1. In
order to have an  acceptable representation of this signal we need
$N=200$ non-zero Fourier coefficients in the non-oversampled case
(corresponding to considering $a=1$). If instead we consider
$a=0.5$, we duplicate the number of coefficients used to represent
the same signal, but it allows us to additionally transmit a
hidden code. In our example we will send a code $\vec{h}$
consisting of $K=12$ numbers, the first 5 digits of which are
shown in the first column of Table I.
\\ \\
{\bf{Case 1:}} If transmitted through a noise-free channel, the
proposed encoding$/$decoding system is, of course, capable 
of transmitting a
great deal of hidden information with accuracy limited only by
machine precision of the calculations involved. The absolute 
value of the coefficients
$\vec{c}$ conveying the information on the chirp signal are
plotted in Figure 2a while the hidden code coefficients $\vec{c'}$
are plotted in Figure 2b. The absolute value 
of the transmitted coefficients $\vec{c''}
= \vec{c} + \vec{c'}$ are those of Figure 2c. For the sake of
avoiding a notorious distortion of coefficients $\vec{c}$ we have
diminished the magnitude of coefficients $\vec{c'}$  by
multiplying by an appropriate scaling factor, which in the absence
of noise can be arbitrarily small. As expected, our decoding
system is capable of reconstructing the hidden code up to the
precision of the numerical calculation, although for space
limitation reasons we
have shown only 5 digits (see the second column of Table I).\\
\\
{\bf{Case 2:}} Now let us illustrate the effect of adding zero
mean random Gaussian noise to the transmitted coefficients. As
would be expected, the quality of the recovery of our signal and
hidden code depends on the variance of the noise ($\sigma^2$)
relative to the size of the signal and hidden code coefficients
respectively. For this case and the following one $\sigma^{2}$ was
fixed such that the signal to noise ratio (in terms of the signal
coefficients) was 40 dB. This allows a reasonably good recovery of
the signal with a small amount of noise distortion. \\
Now that the
noise is fixed relative to the signal, we consider rescaling the
hidden code coefficients $\vec{c'}$ by a constant to achieve the
desired accuracy in the recovery of the code.  We define the
variance ratio as
$$\rho= \frac{\sigma^2}{||\vec{c'}||/{\mbox{Dim}}(\vec{c'})},$$
where ${\mbox {Dim}}(\vec{c'})$  denotes the dimension of vector
$\vec{c'}$. In order to recover 4 significant digits of the code,
we need a variance ratio of about $10^{-5}$ (see third column of
Table I). To achieve such a ratio the magnitude of $\vec{c'}$ has 
to be increased, thereby the absolute value of the 
 transmitted coefficients look as in
Figure 2c. It is important to remark that, although the need to
magnify the coefficients $\vec{c'}$ makes them ``visible'' during
the transmission, this has no effect whatsoever on the signal
reconstruction. Let us recall that $\vec{c'}$ is by definition in
Null($G$), so these coefficients cannot affect the signal in any
way. Hence in this case the hidden code coefficients $\vec{c'}$
actually play a double role. On one hand they cover the
coefficients $\vec{c}$ conveying the information for recovering
the signal $f$ and on the other hand they convey the
information containing the hidden code.\\
\\
{\bf{Case 3:}} Let us finally discuss the limitation of the
proposed system if one does not want the hidden code coefficients
$\vec{c'}$ to dominate the value of the transmitted ones
$\vec{c''}$. With $\sigma^{2}$ fixed as described in the previous
case, the variance ratio was increased to 0.002 (by scaling
$\vec{c'}$ appropriately) and the first digit of the code was
still safely recovered(see the fourth column of Table I). The
absolute value of the 
transmitted coefficients in this case are shown in Figure 2d. 
Note that the transmitted coefficient are dominated now 
by the signal coefficients.

\section{Conclusions}
Oversampling DFT has been considered as providing a means
for transmitting hidden information.
A basic encoding/decoding
system has been discussed. The proposed scheme aims at
transmitting an arbitrary signal
and, simultaneously, embedding a hidden code.
It is important to stress once again  that
the purpose was to discuss in the simplest possible way
the possibility of using redundant transformations
for embedding a hidden code while transmitting an arbitrary
signal. To such an end the idea
has being presented in the context of DFT oversampling.
However, many other redundant transformations could have been
considered by an equivalent treatment. Furthermore, different
ways of embedding a hidden code in the redundant
coefficients could be envisaged.
We feel then confident that
the simple scheme we have introduced here
will stimulate further research in the subject.

\section*{Acknowledgements}
We wish to thank David Lowe for useful discussions and
corrections on the manuscript.\\
Support from EPSRC (GR$/$R86355$/$01) is acknowledged.

\newpage

{\bf {Figure Captions:}}\\

Figure 1: Chirp signal to be transmitted.\\

Figures 2a-2e, from top to bottom.  2a depicts the absolute 
value of the signal
coefficients used in all three cases, 2b the hidden code
coefficients scaled as used in Case 1, and 2c the absolute 
value of the transmitted
coefficients in Case 1. 
Figure 2d represents the absolute value of the 
transmitted coefficients in Case 2 and 2e
the absolute value of the  transmitted coefficients in Case 3.

\newpage

\begin{table}
\begin{center}
\subsection*{}
\begin{tabular}{||c| c c c ||}  \hline
& Case 1 & Case 2 & Case 3\\
code &$\rho= 0$ & $\rho =10^{-5}$& $\rho=0.002$\\
\hline
 3.1492 & 3.1492 & 3.1496 & 3.2286\\
 2.1271 & 2.1271 & 2.1270 & 2.1157\\
 5.1312 & 5.1312 & 5.1316 & 5.2206\\
 1.2835 & 1.2835 & 1.2836 & 1.2939\\
 7.7976 & 7.7976 & 7.7979 & 7.8660\\
 3.7160 & 3.7160 & 3.7164 & 3.7999\\
 8.4139 & 8.4139 & 8.4140 & 8.4360\\
 1.9791 & 1.9791 & 1.9791 & 1.9885\\
 0.5863 & 0.5863 & 0.5868 & 0.6785\\
 5.8321 & 5.8321 & 5.8317 & 5.7570\\
 8.1032 & 8.1032 & 8.1032 & 8.1124\\
 6.4908 & 6.4908 & 6.4907 & 6.4718\\
\hline
\end{tabular}
\end{center}
\vspace{6cm} \caption{Recovered code with given variance ratio.}
\end{table}

\newpage

\begin{figure}

\begin{center}
{\epsfig{figure=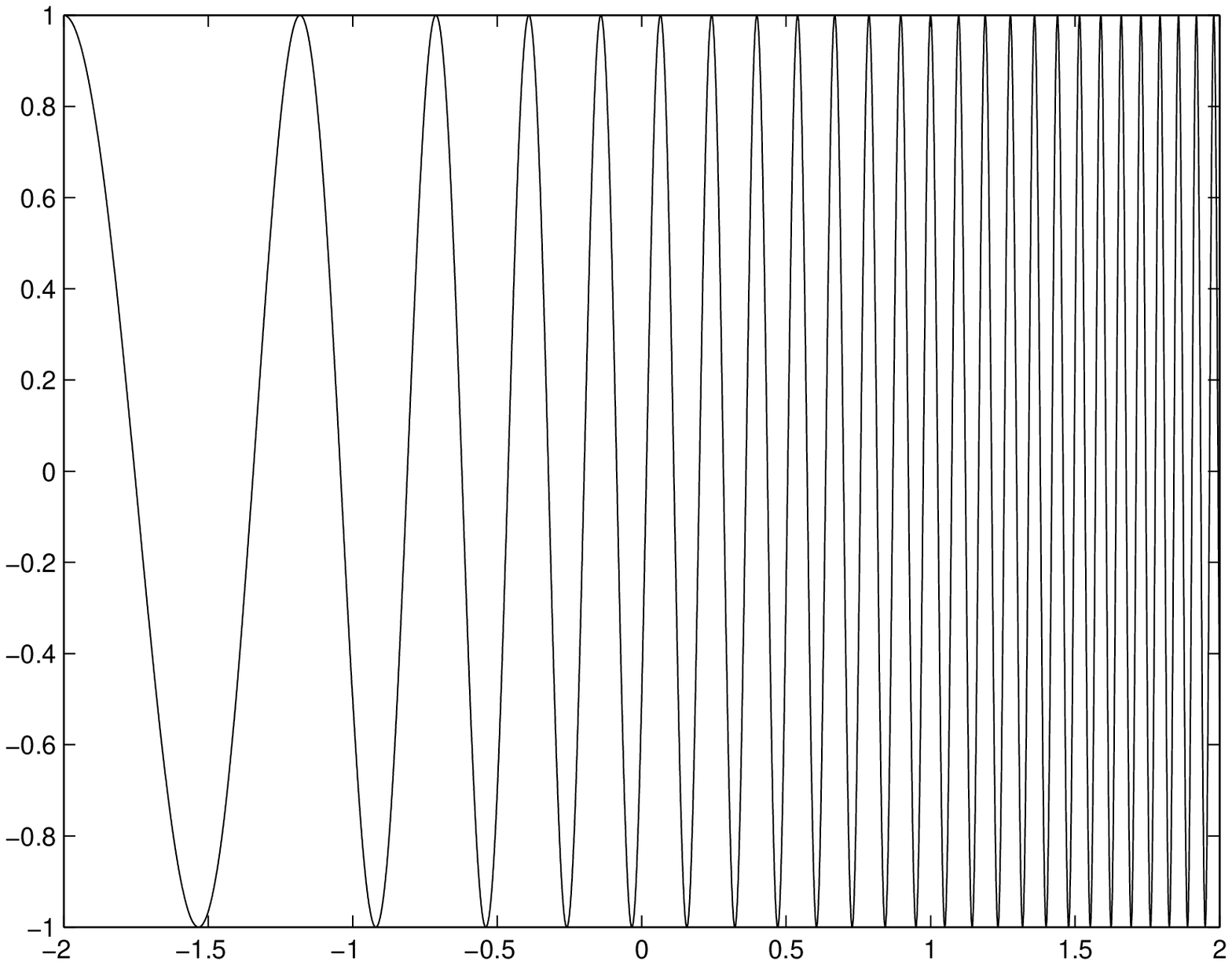,width=.5\textwidth}}
\end{center}
\vspace{6cm} \caption{}
\label{fig1}
\end{figure}

\begin{figure}
\begin{center}
{\epsfig{figure=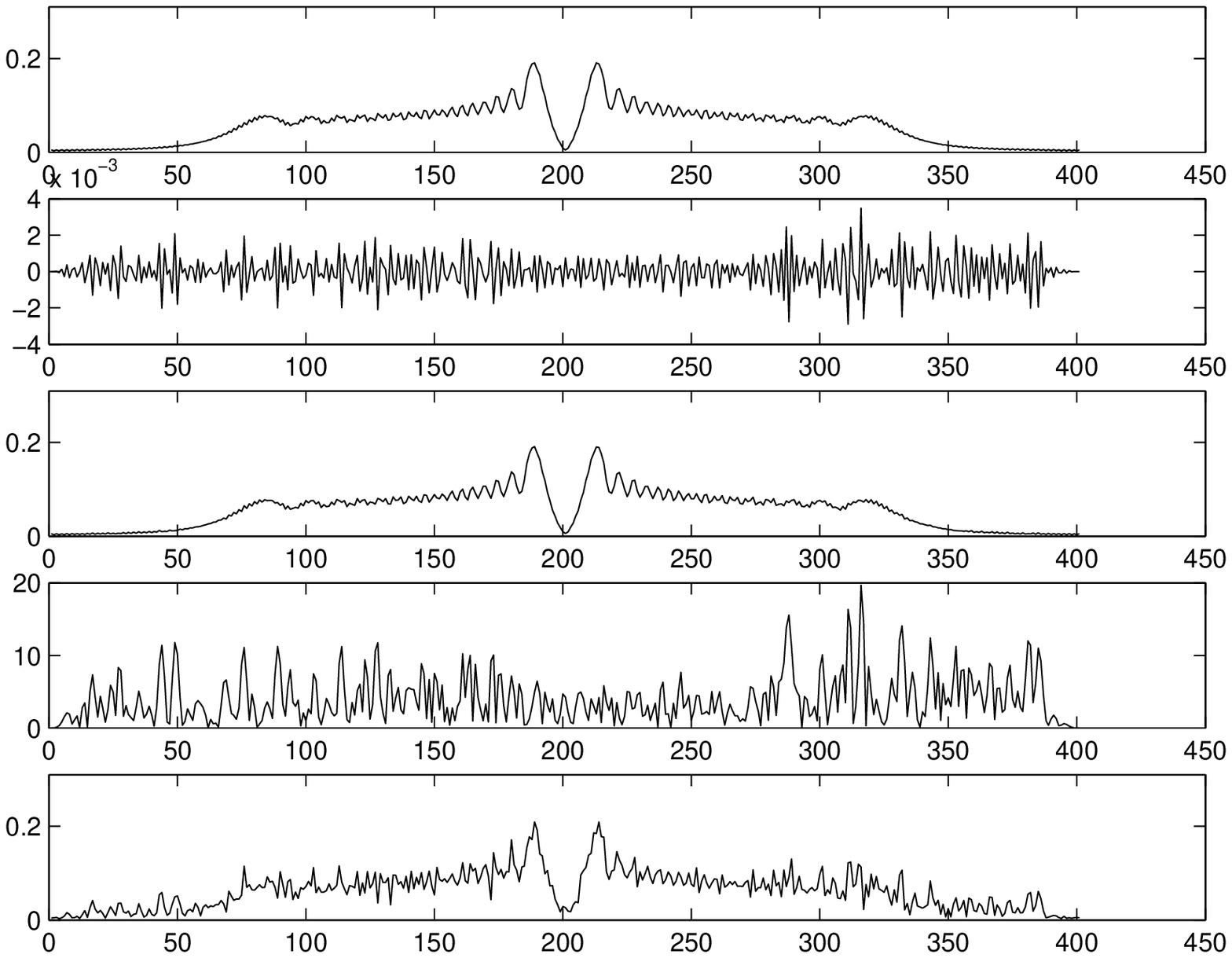,width=.8\textwidth}}
\end{center}
\vspace{6cm} \caption{}
\label{fig2}
\end{figure}

\end{document}